\documentclass[reqno,12pt]{amsart}
\usepackage{amsmath}
\usepackage{amsthm}
\usepackage{amssymb}

\newtheorem{theorem}{Theorem}
\newtheorem{lemma}{Lemma}
\newtheorem{corollary}{Corollary}
\newtheorem{conjecture}{Conjecture}

\begin{document}

\bibliographystyle{amsplain}

\title{Bounds on ternary cyclotomic coefficients}
\author{Bart\l{}omiej Bzd\c{e}ga}
\address{Ul. Kili\'nskiego 19, 62-025 Kostrzyn WLKP, Poland}
\email{\tt exul@wp.pl}
\date{13 December 2008}

\abstract
We present a new bound on $A = \max_n |a_{pqr}(n)|$, where $a_{pqr}(n)$ are the coefficients of a ternary cyclotomic
polynomial $\Phi_{pqr}(x) = \prod_{(k, pqr)=1; \; 0<k<pqr}(x- \xi_{pqr}^k)$ with $p$, $q$, $r$ prime, $p < q,r$, $q \neq r$. We also prove that $|a_{pqr}(n) - a_{pqr}(n-1)| \leqslant 1$.
\endabstract

\maketitle

\section{Introduction}

Let
$$\Phi_{pqr}(x) = \prod_{(k, pqr) = 1; \; 0 < k < pqr}(x- \xi_{pqr}^k) = \sum_n a_{pqr}(n)x^n$$
be a ternary cyclotomic polynomial with $p < q, r$ prime and pairwise different. The coefficients of $\Phi_{pqr}$ have been a subject of studies
for over a century. The main problem was to estimate the following parameters:
\begin{equation} \label{A+-}
A_+ = \max_n a_{pqr}(n); \quad A_- = \min_n a_{pqr}(n); \quad A = \max\{A_+, -A_-\}.
\end{equation}

The first bound on $A$ was given by Bang \cite{bang} who showed that
$$A \leqslant p-1.$$

This bound was improved later by Beiter \cite{beiter}. She proved that
$$A \leqslant p - \Big\lfloor \frac{p}{4} \Big\rfloor.$$

Beiter also came up with a following conjecture:

\begin{conjecture} \label{HBold}
$A \leqslant \frac{p+1}{2}$,
\end{conjecture}

\noindent
now known to be false. Gallot and Moree \cite{gallotmoree}
found infinitely many pairs of primes $q,r$ for every $\varepsilon > 0$ and $p$ sufficiently large, such that
$$A > \left(\frac{2}{3} - \varepsilon\right)p,$$
It updates Beiter's Conjecture into the following form:

\begin{conjecture} \label{HB}
$A \leqslant \frac{2}{3}p$.
\end{conjecture}

This is still an open problem.

In this paper we derive a new bounds on ternary cyclotomic coefficients, which depend on the inverses of $q$ and $r$ modulo $p$ (denoted here by $q'$ and $r'$ respectively). The main result of this paper are the following theorems:

\begin{theorem} \label{M}
Let $A_+$ and $A_-$ be defined as in (\ref{A+-}). Then
$$A_+ \leqslant \min\{2\alpha+\beta,p-\beta\}; \quad -A_- \leqslant \min\{p+ 2\alpha - \beta, \beta\},$$
where $\alpha = \min\{q', r', p-q', p-r'\}$ and $\alpha \beta q r \equiv 1 \; (\text{mod } p)$, $0 < \beta < p$.
\end{theorem}

\begin{theorem} \label{A}
Put $\beta^* = \min\{\beta, p - \beta\}$. Then
\begin{equation} \label{BB}
A \leqslant \min\{2\alpha + \beta^*, p - \beta^*\}.
\end{equation}
\end{theorem}

Theorem \ref{A} improves the bound on $A$ obtained by Bachman \cite{bachman}:
\begin{equation} \label{B}
A \leqslant \min\left\{\frac{p-1}{2}+\alpha, p-\beta^*\right\}.
\end{equation}

One can deduce by \emph{reductio ad absurdum}, that the bound (\ref{BB}) is at least as strong as (\ref{B}). It is also easy to check, that the bound (\ref{BB}) is sharply stronger than (\ref{B}) if and only if $\alpha + \beta^* < \frac{p-1}{2}$, what gives exactly $\frac{1}{2}(p-3)(p-5)$ of all the $(p-1)^2$ pairs $(x,y)$ of residue classes $q$ and $r$ modulo $p$.

As an application, we prove a density result showing that Conjecture \ref{HB} holds for at least $\frac{25}{27}$ of all the ternary cyclotomic polynomials and prove that average $A$ of all the ternary cyclotomic polynomials $\Phi_n$ with the smallest prime factor of $n$ equal to $p$ does not exceed $\frac{p+1}{2}$ (according to Bachman's Theorem these values was $\frac{8}{9}$ and $\frac{7p-1}{12} + o(1)$ respectively).

We also reveal for every prime $p > 3$ some new classes of ternary cyclotomic polynomials $\Phi_{pqr}$ for which the set of coefficients is very small. For example $A \leqslant 3$ if $q \equiv \pm 1 (\text{mod } p)$ and $r \equiv \pm 1 (\text{mod } p)$.

Our method also leads to a simpler proof of the so called \emph{jump one ability} of the ternary cyclotomic coefficients due to Gallot and Moree \cite{gallotmoree2}. It was shown by the present author independently of Gallot and Moree.

\begin{theorem} \label{1}
If $\Phi_{pqr}(x) = \sum_{n \in \mathbb{Z}}a_{pqr}(n)x^n$ is a ternary cyclotomic polynomial, then
$$|a_{pqr}(n) - a_{pqr}(n-1)| \leqslant 1$$
for every $n \in \mathbb{Z}$.
\end{theorem}

\section{The numbers $F_k$}

We define special numbers, which are the key tools in the proof of Theorem \ref{M} and \ref{1}.
Throughout the paper we assume that $k \in \mathbb{Z}$, fix $p$, $q$, $r$ and denote by $a_k$, $b_k$, $c_k$ the unique integers such that $0 \leqslant a_k < p$, $0 \leqslant b_k < q$, $0 \leqslant c_k < r$ and
$$k \equiv a_k qr + b_k rp + c_k pq \; (\text{mod } pqr).$$
Let
$$F_k = \frac{a_k}{p} + \frac{b_k}{q} + \frac{c_k}{r} - \frac{k}{pqr}.$$
Observe that $F_k \in \{0,1,2\}$ for $-(qr+rp+pq) < k < pqr$, since
$$0 \leqslant a_k qr + b_k rp + c_k pq - k < (p-1)qr + (q-1)rp + (r-1)pq + qr + rp + pq = 3pqr.$$
In the remainder of this section we establish the properities of the sequence $F_k$.

\begin{lemma} \label{02}
If $F_k = 0$ then $a_k \leqslant \big\lfloor \frac{k}{qr} \big\rfloor$. If $F_k = 2$ then $a_k \geqslant \big\lceil \frac{k+pq+rp}{qr} \big\rceil$.
\end{lemma}

\begin{proof}
The first implication is obvious. For the second one we note that
$$k+2pqr = a_k qr + b_k rp + c_k pq \leqslant a_k qr + (q-1)rp + (r-1)pq,$$
thus $a_k qr \geqslant k+rp+pq$ and finally $a_k \geqslant \big\lceil \frac{k+pq+rp}{qr} \big\rceil$.
\end{proof}

\begin{lemma} \label{+-}
Let $p_q'$, $p_r'$ be the inverses of $p$ modulo $q$ and $r$ respectively. Then
$$F_k - F_{k-q} = \left\{
\begin{array}{rl}
-1, & \text{if} \quad a_k < r' \text{ and } c_k < p_r';\\
1, & \text{if} \quad a_k \geqslant r' \text{ and } c_k \geqslant p_r';\\
0, & \text{otherwise}.
\end{array}
\right.$$
Analogous statement holds for $F_k - F_{k-r}$ with $c_k$, $r'$, $p_r'$ replaced by $b_k$, $q'$, $p_q'$ respectively.
\end{lemma}

\begin{proof}
Observe that $a_{k-q} \equiv a_k - r' \; (\text{mod } p)$, $c_{k-q} \equiv c_k - p_r' \; (\text{mod } r)$ and $b_{k-q} = b_k$. Therefore
$$a_k - a_{k-q} = \left\{
\begin{array}{rl}
r'-p, & \text{if} \quad a_k < r';\\
r', & \text{if} \quad a_k \geqslant r'.\\
\end{array}
\right.$$
and
$$c_k - c_{k-q} = \left\{
\begin{array}{rl}
p_r'-r, & \text{if} \quad c_k < p_r';\\
p_r', & \text{if} \quad c_k \geqslant p_r'.\\
\end{array}
\right.$$
Let $[P] \in \{0, 1\}$ be the logical value of an expression $P$. Then
\begin{eqnarray*}
F_k - F_{k-q} & = & \frac{a_k - a_{k-q}}{p} + \frac{b_k - b_{k-q}}{q} + \frac{c_k - c_{k-r}}{r} - \frac{1}{pr} \\
& = & \frac{r'}{p} + \frac{p_r'}{r} - \frac{1}{pr} - [a_k < r'] - [c_k < p_r'] \\
& = & 1 - [a_k < r'] - [c_k < p_r'].
\end{eqnarray*}
\end{proof}

\begin{lemma} \label{+--+}
Let $M = \max\{q',r'\}$ and $m = \min \{q',r'\}$. Then
$$F_k - F_{k-q} - F_{k-r} + F_{k-q-r} = \left\{
\begin{array}{rl}
0, & \text{if} \quad a_k < M+m-p;\\
-1, & \text{if} \quad M+m-p \leqslant a_k < m;\\
0, & \text{if} \quad m \leqslant a_k < M;\\
1, & \text{if} \quad M \leqslant a_k < M+m;\\
0, & \text{if} \quad M+m \leqslant a_k.
\end{array}
\right.$$
This Lemma works also for any permutation of $(p, q, r)$ with similarly defined $M$ and $m$.
\end{lemma}

\begin{proof}
Using the method similar to the proof of Lemma \ref{+-}, we obtain
\begin{eqnarray*}
F_k - F_{k-q} - F_{k-r} + F_{k-q-r} & = & - [a_k < q'] - [a_k < r'] + [a_k < q'+r'] + [a_k < q'+r'-p] \\
& = & [a_k < M+m-p] - [a_k < m] - [a_k < M] + [a_k < M+m].
\end{eqnarray*}
Now it is easy to verify the lemma, since $M+m-p < m \leqslant M < M+m$.
\end{proof}

\begin{lemma} \label{8+-}
$$F_k + F_{k-p-q} + F_{k-q-r} + F_{k-r-p} = F_{k-p} + F_{k-q} + F_{k-r} + F_{k-p-q-r}.$$
\end{lemma}

\begin{proof}
By Lemma \ref{+--+}, the value of $F_k - F_{k-q} - F_{k-r} + F_{k-q-r}$ depends only on $k$ modulo $p$. Thus
$$F_k - F_{k-q} - F_{k-r} + F_{k-q-r} = F_{k-p} - F_{k-p-q} - F_{k-r-p} + F_{k-p-q-r}.$$
\end{proof}

\section{Proof of Theorem \ref{M}}

Bloom \cite{bloom} described a relation between the ternary cyclotomic coefficients and the numbers $k$ such that $k = a_k qr + b_k rp + c_k pq$
with $a_k$, $b_k$ and $c_k$ defined in the previous section.
This equality holds if and only if $F_k = 0$, so we can express his result in terms of $F_k$.

\begin{lemma} \label{N}
Denote by $N_d(t_1, t_2, ..., t_l)$ the number of $d$'s in the sequence $(t_1, t_2, ..., t_l)$. Then
\begin{eqnarray*}
a_{pqr}(n) & = & \sum_{k=n-p+1}^{n}(N_0(F_k, F_{k-q-r}) - N_0(F_{k-q}, F_{k-r})) \\
& = & \sum_{k=n-p+1}^{n}(N_2(F_k, F_{k-q-r}) - N_2(F_{k-q}, F_{k-r})) \\
& = & \frac{1}{2} \sum_{k=n-p+1}^{n}(N_1(F_{k-q}, F_{k-r}) - N_1(F_k, F_{k-q-r})).
\end{eqnarray*}
\end{lemma}

\begin{proof}
The first equality is due to Bloom \cite{bloom}. Here we rewrite his proof which uses the formal series:
\begin{eqnarray*}
\Phi_{pqr}(x) & = & \frac{(1-x^{pqr})(1-x^p)(1-x^q)(1-x^r)}{(1-x)(1-x^{qr})(1-x^{rp})(1-x^{pq})} \\
& \equiv & (1-x^q)(1-x^r)(1+x+...+x^{p-1})\sum_{a,b,c \geqslant 0}x^{aqr+brp+cpq} \; (\text{mod } x^{pqr}).
\end{eqnarray*}
Note that if $k \leqslant \deg(\Phi_{pqr}) < pqr$ then there exists at most one triple $(a,b,c)$ such that $k = aqr + brp + cpq$.
This equality holds if and only if $F_k=0$ with $a=a_k$, $b=b_k$, $c=c_k$. Then
\begin{eqnarray*}
a_{pqr}(n) & = & \sum_{k=n-p+1}^n ([F_k=0] - [F_{k-q} = 0] - [F_{k-r}=0] + [F_{k-q-r}=0]) \\
& = & \sum_{k=n-p+1}^{n}(N_0(F_k, F_{k-q-r}) - N_0(F_{k-q}, F_{k-r})).
\end{eqnarray*}
Now in order to simplify the expressions we will use the following notations:
\begin{eqnarray*}
N_0^+ & = & N_0(F_n, F_{n-1}, ..., F_{n-p+1}, F_{n-q-r}, F_{n-q-r-1}, ..., F_{n-q-r-p+1}); \\
N_0^- & = & N_0(F_{n-q}, F_{n-q-1}, ..., F_{n-q-p+1}, F_{n-r}, F_{n-r-1}, ..., F_{n-r-p+1}),
\end{eqnarray*}
and similarly $N_1^+$, $N_1^-$, $N_2^+$, $N_2^-$. We have just proved, that $a_{pqr}(n) = N_0^+ - N_0^-$. Now by Lemma \ref{+--+} we have
\begin{eqnarray*}
N_1^+ + 2N_2^+ - N_1^- - 2N_2^- & = & \sum_{k=n-p+1}^n(F_k - F_{k-q} - F_{k-r} + F_{k-q-r}) \\
& = & \min\{M+m, p\} - M + m - \max\{M+m-p,0\} = 0.
\end{eqnarray*}
Moreover
$$N_0^+ + N_1^+ + N_2^+ = N_0^- + N_1^- + N_2^- = 2p.$$
By simple arithmetical operations, these equalities lead to
$$a_{pqr}(n) = N_0^+ - N_0^- = N_2^+ - N_2^- = \frac{1}{2}(N_1^- - N_1^+).$$
\end{proof}

Using the first equality of Lemma \ref{N}, we consider the 4-tuples $(F_k, F_{k-q}, F_{k-r}, F_{k-q-r})$, where $k \in \{n,n-1,...,n-p+1\}$,
such that $N_0(F_k, F_{k-q-r}) \neq N_0(F_{k-q}, F_{k-r})$. Lemmas \ref{+-} and \ref{+--+} will help us to exclude the existence of
most of the 81 possible 4-tuples.

If $N_0(F_k, F_{k-q}, F_{k-r}, F_{k-q-r}) \in \{0,4\}$ then $N_0(F_n, F_{n-q-r}) = N_0(F_{n-q}, F_{n-r})$, so we are not going to consider these
cases. Also if $N_0(F_k, F_{k-q}, F_{k-r}, F_{k-q-r}) = 2$, then $N_0(F_n, F_{n-q-r}) = N_0(F_{n-q}, F_{n-r})$ or
$|F_k - F_{k-q} - F_{k-r} + F_{k-q-r}| \geqslant 2$, what contradicts Lemma \ref{+--+}, therefore this case also does not need to be considered.

To describe the rest of possibilities we need to observe the following facts:

If $N_0(F_k, F_{k-q}, F_{k-r}, F_{k-q-r}) = 3$ then by Lemma \ref{+-} the only non-zero entry here is equal to 1.

If $N_0(F_k, F_{k-q}, F_{k-r}, F_{k-q-r}) = 1$ then $F_l = 0$ for some $l \in \{k,k-q,k-r,k-q-r\}$. By Lemma \ref{+-} we have $F_{l \pm q} = 1$ and $F_{l \pm r} = 1$, where sign $+$ or $-$ depends on the chosen $l$.

All these cases are described in the table below.

\medskip
\begin{center}
\begin{tabular}{|c|c|c|c|} \hline
Case No. & $(F_k, F_{k-q}, F_{k-r}, F_{k-q-r})$ & $F_k - F_{k-q} - F_{k-r} + F_{k-q-r}$ & $N_0(F_n, F_{n-q-r}) - N_0(F_{n-q}, F_{n-r})$ \\ \hline
$1$ & $(0,0,1,0), (0,1,0,0),$ & $-1$ & $1$ \\
& $(0,1,1,1), (1,1,1,0)$ & & \\ \hline
$2$ & $(0,0,0,1), (1,0,0,0),$ & $1$ & $-1$ \\
& $(1,0,1,1), (1,1,0,1)$ & & \\ \hline
$3$ & $(0,1,1,2), (2,1,1,0)$ & $0$ & $1$ \\ \hline
$4$ & $(1,0,2,1), (1,2,0,1)$ & $0$ & $-1$ \\ \hline
\end{tabular}
\end{center}
\medskip

Denote by $C_l$ the number of integers $k \in \{n,n-1,...,n-p+1\}$ for which the $l$th case occurs. Then we have
$$A_+ \leqslant C_1 + C_3; \quad -A_- \leqslant C_2 + C_4.$$
In order to prove Theorem \ref{M} it is enough to show that
\begin{equation} \label{C}
C_1, C_2 \leqslant \alpha; \quad C_3 \leqslant \min\{\alpha + \beta, p-\alpha - \beta\}; \quad C_4 \leqslant \min\{\beta - \alpha, p+\alpha-\beta\}.
\end{equation}

In fact, we will count values of $a_k$ instead of $k$ (there is a bijection between the sets
$\{n,n-1,...,n-p+1\}$ and $\{a_n,a_{n-1},...,a_{n-p+1}\}$, because $a_kqr \equiv k \; (\text{mod } p)$).

Note that $\alpha = \min\{m, p-M\}$, where $M$ and $m$ are defined in Lemma \ref{+--+}.

\medskip \noindent
\textbf{Case 1} \\
By lemma \ref{+--+} we have $M+m-p \leqslant a_k < m$, so
$$C_1 \leqslant m - \max\{0,M+m-p\} = \min\{m,p-M\} = \alpha.$$

\medskip \noindent
\textbf{Case 2} \\
By Lemma \ref{+--+}, we have $M \leqslant a_k < M+m$, so
$$C_2 \leqslant \min\{M+m, p\} - M = \min\{m, p-M\} = \alpha.$$

Note that
$$\text{if} \quad M + m \geqslant p \quad \text{then} \quad \alpha = p - M \text{ and } \beta = p - m$$
and
$$\text{if} \quad M + m \leqslant p \quad \text{then} \quad \alpha = m \text{ and } \beta = M.$$
We also put $\gamma = \Big\lfloor \frac{n}{qr} \Big\rfloor + 1$ and remind that here $k \in \{n,n-1,...,n-p+1\}$.

In order to simplify the notation, we divide the third case into cases $3a$ and $3b$ and define $C_{3a}$ and $C_{3b}$ as above for the 4-tuples $(0,1,1,2)$ and $(2,1,1,0)$ respectively. Obviously, $C_3 = C_{3a} + C_{3b}$.

\medskip \noindent
\textbf{Case 3a} \\
By Lemma \ref{+-} we have here $a_k < m,M$, thus by Lemma \ref{+--+} $a_k < M+m-p$. By Lemma \ref{02}
$$a_k < \gamma$$
and
$$a_k-M-m+2p = a_{k-q-r} \geqslant \gamma.$$
Finally
$$\max\{\gamma + M + m - 2p, 0\} \leqslant a_k < \min\{\gamma, M+m-p\},$$
and we obtain
\begin{eqnarray*}
C_{3a} & \leqslant & \min\{\gamma, M+m-p\} - \max\{\gamma + M + m - 2p, 0\} \\
& = & \min\{\gamma, p-\gamma, M+m-p, 2p-M-m\} \\
& \leqslant & \min\{M+m-p, 2p-M-m\} \\
& = & \min\{\alpha + \beta, p-\alpha - \beta\},
\end{eqnarray*}
as long as $M+m \geqslant p$. Otherwise $C_{3a} = 0$.

\medskip \noindent
\textbf{Case 3b} \\
By Lemma \ref{+-} $a_k \geqslant m,M$, so by Lemma \ref{+--+} $a_k \geqslant M+m$. By Lemma \ref{02}
$$a_k-M-m = a_{k-q-r} < \gamma$$
and
$$a_k \geqslant \gamma.$$
Finally
$$\max\{\gamma, M+m\} \leqslant a_k < \min\{p, \gamma + M + m\}.$$
Therefore
\begin{eqnarray*}
C_{3b} & \leqslant & \min\{p, \gamma+M+m\} - \max\{\gamma, M+m\} \\
& = & \min\{\gamma, p-\gamma, M+m, p-M-m\} \\
& \leqslant & \min\{M+m, p-M-m\} \\
& = & \min\{\alpha + \beta, p - \alpha - \beta\},
\end{eqnarray*}
as long as $M+m \leqslant p$. Otherwise $C_{3b} = 0$.

\medskip \noindent
\textbf{Case 3} \\
Note that cases 3a and 3b are excluding each other. Thus $C_3 \leqslant \min\{\alpha + \beta, p - \alpha - \beta\}$.

\medskip \noindent
\textbf{Case 4} \\
Let us assume that $q'=m$ and $r'=M$. By Lemma \ref{+-}, we have $m \leqslant a_k < M$ (for $F_{k-q}=0$) or $M \leqslant a_k < m$ (when $F_{k-r}=0$).
The second inequality is impossible, so
$$F_{k-q} = 0 \quad \text{and} \quad F_{k-r} =2.$$
By Lemma \ref{02}
$$a_k-m = a_{k-q} < \gamma$$
and
$$a_k-M+p = a_{k-r} \geqslant \gamma.$$
Finally
$$\max\{M+\gamma-p,m\} \leqslant a_k < \min\{m+\gamma, M\},$$
and
\begin{eqnarray*}
C_4 & \leqslant & \min\{m+\gamma, M\} - \max\{M+\gamma-p,m\} \\
& = & \min\{\gamma, p-\gamma, p-M+m, M-m\} \\
& \leqslant & \min\{p-M+m,M-m\} \\
& = & \min\{\beta - \alpha, p+\alpha-\beta\}.
\end{eqnarray*}

That completes the verification of (\ref{C}) and the proof of Theorem \ref{M}. \qed

\section{The bound on $A$}

In this section we derive a bound on $A = \max\{A_+, -A_-\}$. We also establish some infinite families of triples $(p,q,r)$ with restrictions on $q$ and $r$ modulo $p$ only, for which $A$ is bounded by a constant independent of $p$, $q$, $r$.

We also apply our bound on $A$ to estimate the density of the set of ternary cyclotomic polynomials such that $\frac{A}{p} \leqslant c$, for any real $c > 0$. In view of Conjecture \ref{HB}, the most interesting case is $c = \frac{2}{3}$.

At the end we prove a weaker version of the old Beiter's Conjecture.

\begin{proof}[Proof of Theorem \ref{A}]
By Theorem \ref{M}, we have
$$A \leqslant \max\{\min\{2\alpha+\beta, p-\beta\}, \min\{p+2\alpha-\beta, \beta\}\}.$$
If $\beta < \frac{1}{2}p$ then
$$A \leqslant \max\{\min\{2\alpha+\beta, p-\beta\}, \beta\} = \min\{2\alpha+\beta, p-\beta\} = \min\{2\alpha+\beta^*, p-\beta^*\}.$$
Also if $\beta > \frac{1}{2}p$ then
$$A \leqslant \max\{p-\beta, \min\{p+2\alpha-\beta, \beta\}\} = \min\{2\alpha+p-\beta, \beta\} = \min\{2\alpha+\beta^*, p-\beta^*\}.$$
\end{proof}

\begin{corollary} \label{s}
Let $p > 3$ and $p = 2a+1 = 3b \pm 1 = 4c \pm 1 = 6d \pm 1$ for some integers $a$, $b$, $c$, $d$. If $q$ is congruent to one of the numbers $\pm 1$, $\pm a$, $\pm b$, $\pm c$, $\pm d$ modulo $p$ and also $r$ is congruent to one of them modulo $p$, then $A \leqslant 18$.
\end{corollary}

\begin{proof}
Observe that in all these cases $\alpha \leqslant \beta^* \leqslant 6$. Then by Theorem \ref{A}, $A \leqslant 2\alpha + \beta^* \leqslant 18$.
\end{proof}

Note that if both $q$ and $r$ are congruent to $\pm 1$ modulo $p$, then $\alpha = \beta^* = 1$ and $A \leqslant 3$.

\begin{corollary} \label{d}
Let $c > 0$ be a real number. Denote by $D(c)$ the density of ternary cyclotomic polynomials for which $\frac{A}{p} < c$. Then
$$D(c) \geqslant \left\{
\begin{array}{rl}
\frac{4}{3}c^2, & \text{if} \quad 0 < c < \frac{1}{2};\\
1 - \frac{2}{3}(3-4c)^2, & \text{if} \quad \frac{1}{2} \leqslant c < \frac{3}{4}.
\end{array}
\right.$$
Also $D(c) = 1$ if $c \geqslant \frac{3}{4}$.
\end{corollary}

\begin{proof}
Let us denote by $P_{p,n}(x,y)$ the probability that $\alpha=x$ and $\beta^*=y$, where $q \neq r$ are random primes from the set $\{p+1, p+2,...,n\}$ and $\alpha$ and $\beta^*$ are computed for the polynomial $\Phi_{pqr}$. Dirichlet's Theorem says that the densities of primes in the arithmetical progressions $1, p+1, 2p+1,...$; $2, p+2, 2p+2,...$; ...; $p-1, 2p-1,3p-1,...$ are the same. Then probabilities
$$P_p(x,y) = \lim_{n \rightarrow \infty} P_{p,n}(x,y)$$
are equal for every integers $1 \leqslant x < y \leqslant \frac{p-1}{2}$. It implies that if $p \rightarrow \infty$ then the distribution of ($\frac{x}{p},\frac{y}{p})$ converges to the uniform distribution over the triangle $T$ on the vertices $(0,0)$, $\left(0,\frac{1}{2}\right)$, $\left(\frac{1}{2}, \frac{1}{2}\right)$.

Note that for random $m \in \{1,2,...,n\}$, $m = pqr$, $p < q,r$, the expected value $\mathbb{E}_n(p) \rightarrow \infty$ when $n \rightarrow \infty$. Then $D(c)$ is not smaller than the area of some polygon divided by the area of the triangle $T$. Precisely, $D(c) \geqslant 8S(c)$,
where $S(c)$ is the area of the polygon defined by inequalities:
$$0 < x < \frac{1}{2}; \quad 0 < y < \frac{1}{2}; \quad x < y; \quad c > \min\{2x+y, 1-y\}.$$
The last inequality is due to Theorem \ref{A}. We can compute $S(c)$ by simple summing the areas of some triangles. We obtain that
$$S(c) = \left\{
\begin{array}{rl}
\frac{1}{6}c^2, & \text{if} \quad 0 < c < \frac{1}{2};\\
\frac{1}{8} - \frac{1}{12}(3-4c)^2, & \text{if} \quad \frac{1}{2} \leqslant c < \frac{3}{4};\\
\frac{1}{8}, & \text{if} \quad c \geqslant \frac{3}{4}.
\end{array}
\right.$$
This completes the proof of Corollary \ref{d}.
\end{proof}

We can apply our estimation of $D(c)$ to check that Conjecture \ref{HB} is true for the set of ternary cyclotomic polynomials of density $\geqslant \frac{25}{27}$. The old Beiter's Conjecture \ref{HBold} holds for at least $\frac{1}{3}$ of all the ternary cyclotomic polynomials.

Although Conjecture \ref{HBold} does not hold in general, we are able to prove a weaker version of it, with the same bound. Let $\overline{A}(p)$ denotes the average value of $A$ of all the ternary cyclotomic polynomials $\Phi_n$ with the smallest prime dividing $n$ equal to $p$.

\begin{corollary} \label{av}
$\overline{A}(p) \leqslant \frac{p+1}{2}$
\end{corollary}

\begin{proof}
Let $a(i,j) = \min\{2\alpha_{i,j} + \beta_{i,j}^*, p - \beta_{i,j}^*\}$, where $\alpha_{i,j}$ and $\beta_{i,j}^*$ are equal to $\alpha$ and $\beta^*$ computed for the polynomial $\Phi_{pqr}$ with $q' \equiv i \; (\text{mod } p)$ and $r' \equiv i \; (\text{mod } p)$. Based on Theorem \ref{A}, using Dirichlet's Theorem as in the proof of Corollary \ref{d}, we obtain
$$ \overline{A}(p) \leqslant \frac{1}{(p-1)^2} \sum_{i=1}^{p-1} \sum_{j=1}^{p-1} a(i,j) = \frac{4}{(p-1)^2} \sum_{i=1}^{(p-1)/2} \sum_{j=1}^{(p-1)/2} a(i,j).$$
Let $k \leqslant \frac{p-1}{2}$ be a nonnegative integer. Then
\begin{eqnarray*}
\sum_{i=1}^k a\left(i, i+ \frac{p-1}{2}-k\right) & = & \sum_{i=k}^l \min\left\{3i-k+\frac{p-1}{2}, k-i+\frac{p+1}{2}\right\} \\
& = & \frac{(p+1)k}{2} + \sum_{i=1}^k \min\{3i-k-1, k-i\} \\
& = & \frac{(p+1)k}{2}.
\end{eqnarray*}
It implies that $\overline{A}(p) \leqslant \frac{p+1}{2}$.
\end{proof}

\section{Proof of Theorem \ref{1}}

First we present a simple expression for the difference of the two consecutive coefficients of a ternary cyclotomic polynomial in terms of $F_k$:

\begin{lemma} \label{n-1}
Put
$$N_+ = N_1(F_n,F_{n-p-q},F_{n-q-r},F_{n-r-p})$$
and
$$N_- = N_1(F_{n-p},F_{n-q},F_{n-r},F_{n-p-q-r}).$$
Then
$$a_{pqr}(n) - a_{pqr}(n-1) = \frac{1}{2}(N_- - N_+).$$
Moreover
\begin{eqnarray*}
a_{pqr}(n) - a_{pqr}(n-1) & = & N_0(F_n,F_{n-p-q},F_{n-q-r},F_{n-r-p}) - N_0(F_{n-p},F_{n-q},F_{n-r},F_{n-p-q-r}) \\
& = & N_2(F_n,F_{n-p-q},F_{n-q-r},F_{n-r-p}) - N_2(F_{n-p},F_{n-q},F_{n-r},F_{n-p-q-r}).
\end{eqnarray*}
\end{lemma}

\begin{proof}
By Lemma \ref{N}
\begin{eqnarray*}
a_{pqr}(n) - a_{pqr}(n-1) & = & \frac{1}{2} \sum_{k=n-p+1}^{n}(N_1(F_{k-q}, F_{k-r}) - N_1(F_k, F_{k-q-r})) \\
& & - \frac{1}{2} \sum_{k=n-p}^{n-1}(N_1(F_{k-q}, F_{k-r}) - N_1(F_k, F_{k-q-r})) \\
& = & \frac{1}{2}(N_1(F_{n-p},F_{n-q},F_{n-r},F_{n-p-q-r}) - N_1(F_n,F_{n-p-q},F_{n-q-r},F_{n-r-p})) \\
& = & \frac{1}{2}(N_- - N_+).
\end{eqnarray*}
The remaining two equalities can be shown in the same way.
\end{proof}

Now we are ready to prove Theorem \ref{1}. By Lemma \ref{n-1} we have
$$|a_{pqr}(n) - a_{pqr}(n-1)| = \frac{1}{2}|N_- - N_+| \leqslant 2,$$
where equality may hold only if $N_- = 4$, $N_+ = 0$ or $N_+ = 4$, $N_- = 0$. We will show that it is impossible.

Indeed, for some permutation $(t,u,v)$ of $(p,q,r)$ by Lemma \ref{n-1} we have $F_{n-t} = F_{n-u} \in \{0,2\}$ in case of
$(F_n, F_{n-p-q}, F_{n-q-r}, F_{n-r-p}) = (1, 1, 1, 1)$. Therefore $|F_n - F_{n-t} - F_{n-u} + F_{n-t-u}| = 2$.
Also if $(F_{n-p}, F_{n-q}, F_{n-r}, F_{n-p-q-r}) = (1, 1, 1, 1)$ then for some permutation $(t,u,v)$ we have $F_{n-t-u} = F_{n-u-v} \in \{0,2\}$ and
$|F_{n-u} - F_{n-t-u} - F_{n-u-v} - F_{n-t-u-v}| = 2$. Both cases contradict Lemma \ref{+--+}. This completes the proof of Theorem \ref{1}. \qed

\section*{Acknowledgments}
The research was done when the author was a student at the Faculty of Mathematics and Computer Science on the Adam Mickiewicz University in Pozna\'n. The author would like to thank Wojciech Gajda for suggesting the problem and his help in improving the paper. He also would like to thank Pieter Moree for helpful comments.

\providecommand{\bysame}{\leavevmode\hbox to3em{\hrulefill}\thinspace}
\providecommand{\MR}{\relax\ifhmode\unskip\space\fi MR }
\providecommand{\MRhref}[2]{\href{http://www.ams.org/mathscinet-getitem?mr=#1}{#2}}
\providecommand{\href}[2]{#2}

\end{document}